\documentclass[reqno,12pt]{amsart}

\usepackage{epsf}
\usepackage{graphics}
\usepackage{amssymb}
\usepackage{amsmath}

\date{}

\theoremstyle{plain}
\newtheorem{theorem}{Theorem}

\theoremstyle{definition}

\theoremstyle{remark}

\newtheorem*{remark}{Remark}

\def\N{{\mathbb N}}

\title{Gordian distance and Vassiliev invariants} 

\author{Sebastian Baader}

\begin{document}

\begin{abstract} The Gordian distance between two knots measures how many crossing changes are needed to transform one knot into the other. It is known that there are always infinitely many non-equivalent knots `between' a pair of knots of Gordian distance two. In this paper we prove an extreme generalisation of this fact: there are knots with arbitrarily prescribed Vassiliev invariants between every pair of knots of Gordian distance two.
\end{abstract}

\maketitle

\section{Introduction}

A crossing change is an operation on knots that allows us to untie every knot in finitely many steps. We say that two knots differ by a crossing change, if they are related by a strand passage operation along an embedded disc with one singularity, as illustrated in figure~1. The Gordian distance $d_G$ between two knots is the least number of crossing changes needed to transform one knot into the other. 

\begin{figure}[ht]
\scalebox{1}{\raisebox{-0pt}{$\vcenter{\hbox{\epsffile{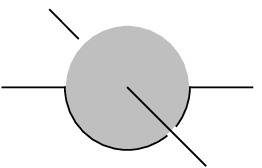}}}$}} \quad
$\longleftrightarrow$ \quad
\scalebox{1}{\raisebox{-0pt}{$\vcenter{\hbox{\epsffile{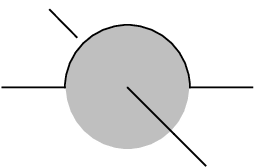}}}$}}
\caption{}
\end{figure}

The Gordian distance gives rise to a simplicial complex, the \emph{Gordian complex} of knots: every isotopy class of 
smooth oriented knots in $S^3$ corresponds to a vertex. Two vertices are connected by an edge, if the corresponding knots have Gordian distance one. Similarly, an $n$-simplex of the Gordian complex is a set of $n+1$ knots whose pairwise Gordian distance is one. The structure of the Gordian complex is rather complicated. As Gambaudo and Ghys proved in \cite{GG}, it contains quasi-isometric embeddings of high rank grids. Moreover, every edge of the Gordian complex is contained in a simplex of infinite dimension (see~\cite{HU}).

Given two knots of Gordian distance two, we may ask how many knots are `between' them. Quite surprisingly, the answer is `infinitely many', as was recently shown in \cite{B}. This can be strengthened a lot.

\begin{theorem} Let $m$ be a natural number and $K$ any knot. For every pair of knots $K_0$ and $K_1$ of Gordian
distance two, there exists a knot $K'$ with the following properties:
\begin{enumerate}
\item the Gordian distance between $K'$ and $K_0$, $K_1$ is one,\\
\item all Vassiliev invariants of $K'$ of order
less than or equal to $m$ coincide with those of $K$.\\
\end{enumerate}
\end{theorem}

Theorem~1 also provides a generalisation of Ohyama's result (\cite{O}), who showed that every finite set of Vassiliev invariants can be realised by an unknotting number one knot. In the next section, we prove the special case $m=2$ of theorem~1; the general case $(m \geq 2)$ is proved in the third section.

\section{The Casson invariant of Knots}

The second coefficient $a_2$ of the Conway polynomial is called the \emph{Casson invariant} of knots. It is the unique Vassiliev invariant of order two on knots, up to multiplication by a constant. The following skein relation provides an effective way of computing $a_2$:
$$a_2(\scalebox{0.2}{\raisebox{8pt}{$\vcenter{\hbox{\epsffile{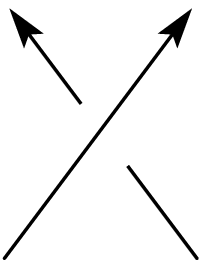}}}$}})
-a_2(\scalebox{0.2}{\raisebox{8pt}{$\vcenter{\hbox{\epsffile{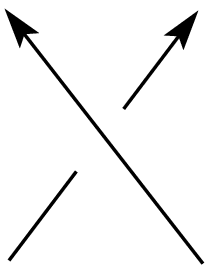}}}$}})
=lk(\scalebox{0.2}{\raisebox{8pt}{$\vcenter{\hbox{\epsffile{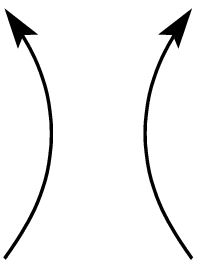}}}$}}).$$
Here $lk(L)$ denotes the linking number of a two component link $L$. 

\begin{proof}[Proof of theorem~1, for $m=2$]
Let $K_0$, $K_1$ be two knots of Gordian distance two. By a standard argument, $K_0$ has a diagram with two neighbouring crossings $A$ and $B$, as shown in figure~2, such that a simultaneous crossing change at $A$ and $B$ transforms $K_0$ into $K_1$. For the moment, we assume that the crossing $A$ is positive. Changing the section of $K_0$, as shown in figure~3, with a certain number $k$ of half-twists above and below the crossing $B$, we define a new knot $K'$, which is `between' $K_0$ and $K_1$. Indeed, a crossing change at $C$ transforms $K'$ into $K_0$, whereas a crossing change at $B$ transforms $K'$ into $K_1$. We note that $C$ is a negative crossing, since $A$ was a positive crossing.

\begin{figure}[ht]
\scalebox{1}{\raisebox{-0pt}{$\vcenter{\hbox{\epsffile{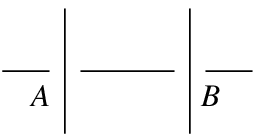}}}$}}
\caption{}
\end{figure}

\begin{figure}[ht]
\scalebox{1}{\raisebox{-0pt}{$\vcenter{\hbox{\epsffile{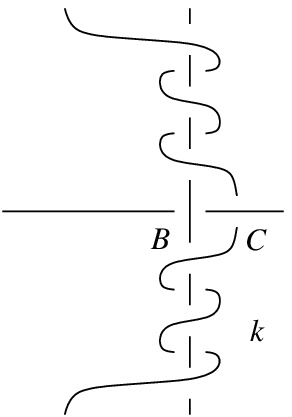}}}$}}
\caption{}
\end{figure}

Applying the skein relation for the Casson invariant at the crossing $C$ of $K'$, we get 
$$a_2(K_0)-a_2(K')=lk(L),$$
where $L$ is the two component link shown in figure~4 (the smoothing must be like this, since the crossing $C$ is negative). Varying the number $k$ of positive (or negative) half-twists, we can arrange
$$lk(L)=a_2(K_0)-N,$$ 
for any given integer $N$. (The two `slalom strands' belong to different components of $L$, whence the vertical strand links with precisely one of the two `slalom strands'. Thus we can vary the linking number of $L$, as we like.) In case $A$ is a negative crossing, we may apply an analogous construction. This proves the special case $m=2$ of theorem~1.

\begin{figure}[ht]
\scalebox{1}{\raisebox{-0pt}{$\vcenter{\hbox{\epsffile{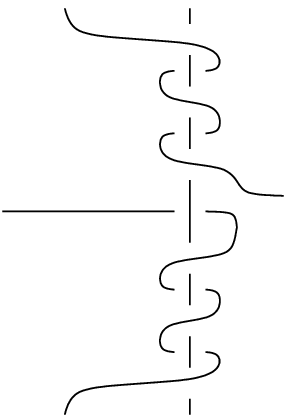}}}$}}
\caption{}
\end{figure}
\end{proof}

\section{$C_n$-moves and Vassiliev invariants}

The proof of Theorem~1 involves Vassiliev invariants and certain $C_n$-moves, 
which were defined by K.~Habiro in~\cite{Ha1} (see also~\cite{Ha2}). 
A $C_1$-move is simply a crossing change. For $n \geq 3$, 
there exist $C_n$-moves and \emph{special} $C_n$-moves. We will only use the latter. 
A special $C_n$-move is defined diagrammatically in figure~5 
(the small loop with an arrow should be ignored at this point).
It takes place in a section with $2(n+1)$ endpoints or $n+1$ strands,
respectively. The strands are numbered from $1$ to $n+1$ and are all connected
outside the indicated section, since they belong to one knot $K$.
Going along $K$ according to its orientation, starting at the first strand, we
encounter the other strands in a certain order which depends on how the
strands are connected outside the indicated section. This order defines a
permutation, say $\sigma \in S_{n}$, of the numbers 2, 3,$\ldots$, n+1.

\begin{figure}[ht]
\scalebox{0.8}{\raisebox{+0pt}{$\vcenter{\hbox{\epsffile{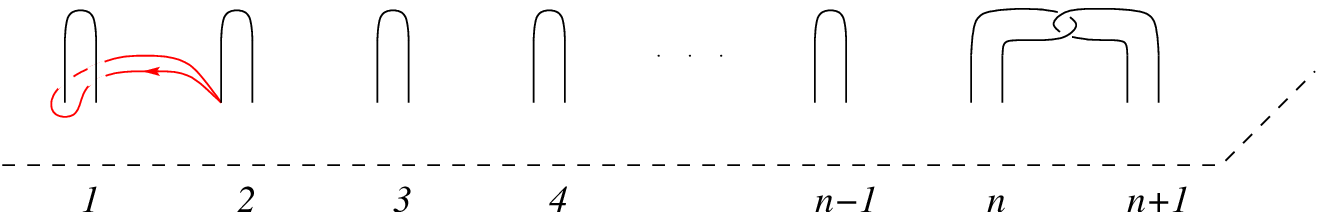}}}$}}

\bigskip
$\updownarrow$ 

\bigskip
\scalebox{0.8}{\raisebox{+0pt}{$\vcenter{\hbox{\epsffile{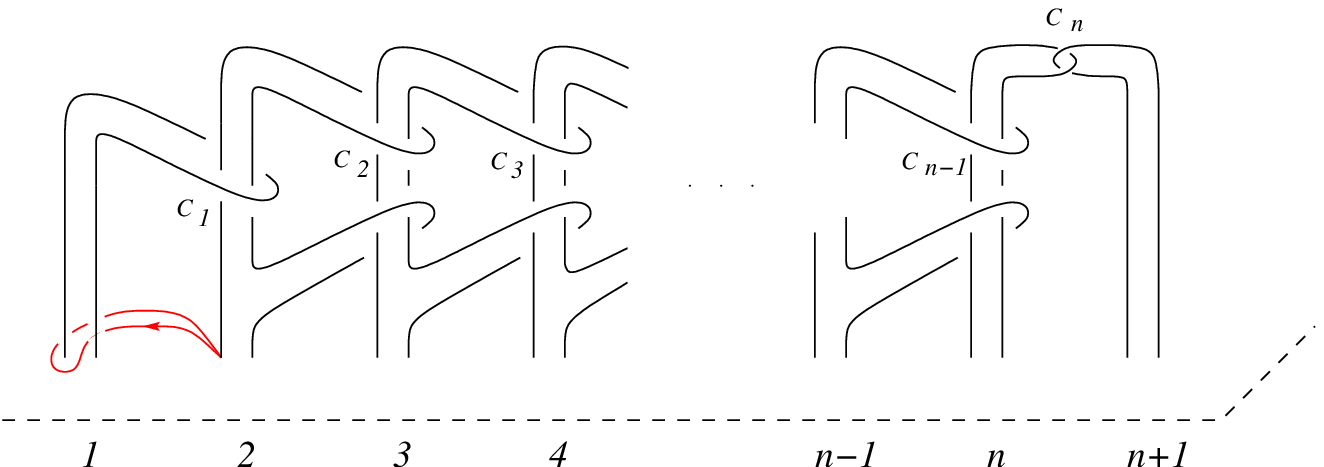}}}$}}
\caption{}
\end{figure}

In~\cite{OTs}, Y.~Ohyama and T.~Tsukamoto explain the effect of a $C_n$-move on
Vassiliev invariants of order $n$. Their result (\cite{OTs}, theorem~1.2)
implies the following:

\begin{enumerate}
\item A $C_n$-move does not change the values of Vassiliev invariants of order
less than $n$.\\
\item Let $K$ and $\widetilde K$ be two knots which differ by one $C_n$-move,
and $v_n$ any Vassiliev invariant of order $n$. Then $v_n(K)-v_n(\widetilde
K)$ depends only on the permutation $\sigma \in S_{n}$ defined by the cyclic
order of the $(n+1)$ strands of the section where the $C_n$-move takes place,
and on the product of the signs of the crossings $c_1$, $c_2$, $c_3$,$\ldots$, $c_n$ (see figure~5). 
\end{enumerate}

The last statement of (2) follows from the proof of theorem~1.2 in~\cite{OTs}.

We will use a slightly modified version of a special $C_n$-move, as shown in figure~6. It is actually equivalent to a special $C_n$-move, as we can see by dragging one end of the second strand in figure~5 along the small loop, in direction indicated by the arrow. This operation relates the two diagrams of figure~5 with those of figure~6.

\begin{figure}[ht]
\scalebox{0.8}{\raisebox{+0pt}{$\vcenter{\hbox{\epsffile{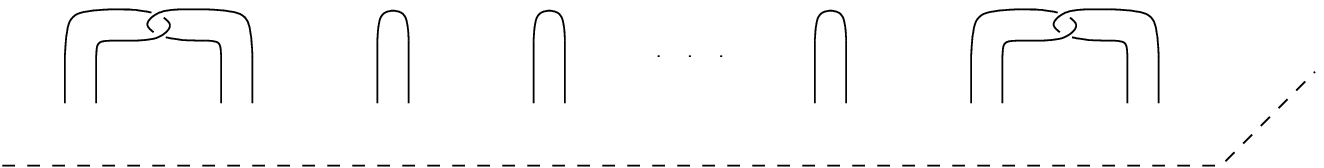}}}$}}

\bigskip
$\updownarrow$ 

\bigskip
\scalebox{0.8}{\raisebox{+0pt}{$\vcenter{\hbox{\epsffile{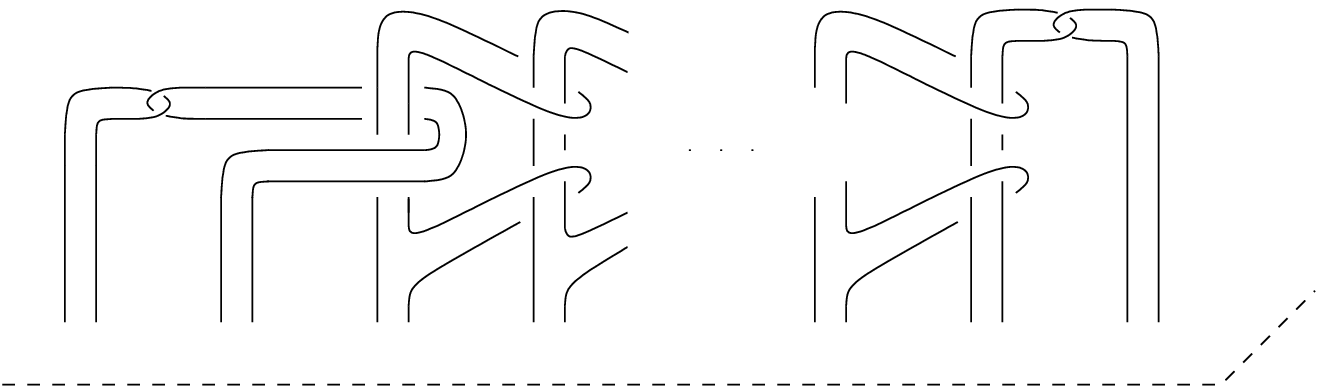}}}$}}
\caption{}
\end{figure}

\begin{proof}[Proof of theorem~1]
Let $K_0$, $K_1$ be two fixed knots of Gordian distance two, as before. By the result of the preceding section, the case $m=2$ of theorem~1 is already settled: 
there exists a knot $K'$ with any prescribed Casson invariant and $d_G(K',K_0)=d_G(K',K_1)=1$. Now let us assume that theorem~1 is true for $m \in \N$, $m \geq 2$, i.e. for any given knot $K$, there exists a knot $K'$, such that:
\begin{enumerate}
\item $d_G(K',K_0)=d_G(K',K_1)=1$,\\
\item all Vassiliev invariants of $K'$ of order less than or equal to $m$ coincide with those of $K$.\\
\end{enumerate}
We will build a new knot $K''$ that meets condition~(2), for $m+1$.

Again, by a standard argument, $K'$ has a diagram with a section as shown in figure~7, with $2$ distinguished spots $A_1$ and $A_2$, where crossing changes should take place to obtain $K_0$ and $K_1$, respectively.
In this section, we add $m-2$ single arcs between the two spots $A_1$ and $A_2$, in an arbitrary way, so that we obtain a section where we can apply a special $C_{m+1}$-move, as shown at the top of figure~6. 

\begin{figure}[ht]
\scalebox{1}{\raisebox{+0pt}{$\vcenter{\hbox{\epsffile{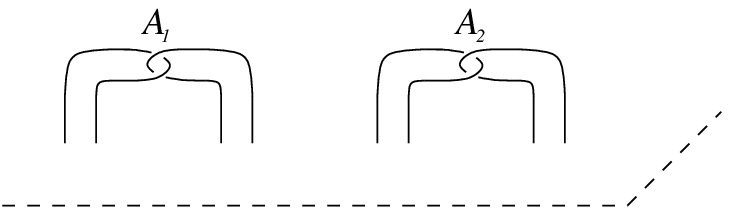}}}$}}
\caption{}
\end{figure}

Going along the knot $K'$ according to its orientation, starting at the leftmost strand,
we encounter the other strands in a certain order, which again defines a
permutation $\sigma \in S_{m+1}$. We claim that we can easily change the order of two strands. 
This is obvious, except for the case of two strands $s_1$, $s_2$ belonging to two different clasps, as
shown at the top of figure~8 (the middle arcs are omitted here). By applying regular isotopy only, we can transform the
indicated section into the section at the bottom of the figure. We observe that the order of
the strands $s_1$, $s_2$ is reversed, there. Moreover, a crossing change at the
clasps involving $s_1$ or $s_2$ yields the same knots as before.

\begin{figure}[ht]
\scalebox{1}{\raisebox{-0pt}{$\vcenter{\hbox{\epsffile{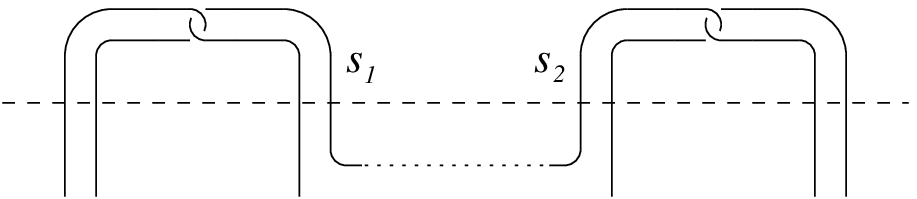}}}$}}

\bigskip
\bigskip
\scalebox{1}{\raisebox{-0pt}{$\vcenter{\hbox{\epsffile{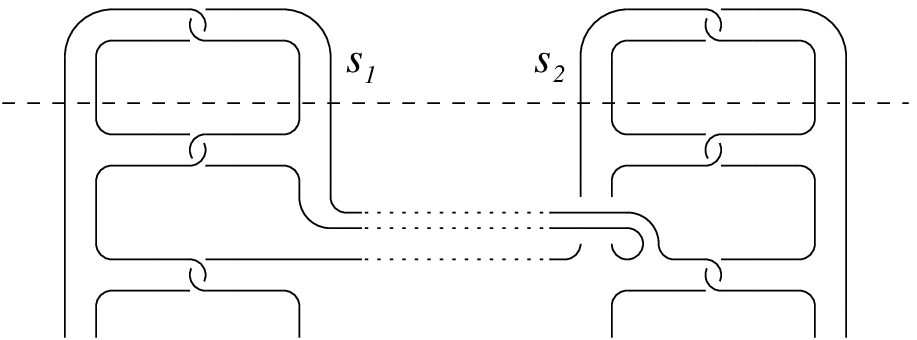}}}$}}
\caption{}
\end{figure}

Hence we may assume that we encounter the strands of figure~6 in an
order corresponding to a fixed prescribed permutation $\sigma \in S_{m+1}$.

Applying the special $C_{m+1}$-move of figure~6, we obtain a knot $\widetilde{K}'$ which is still a neighbour of $K_0$ and $K_1$ in the Gordian complex of knots: $d_G(\widetilde{K}',K_0)=d_G(\widetilde{K}',K_1)=1$.
In order to arrange the Vassiliev invariants of order $\leq m+1$ of $\widetilde{K}'$, we must invoke Habiro's theorem.

\begin{theorem}[Habiro \cite{Ha2}] Two knots have the same Vassiliev invariants up to order $n$, if, and only if, they are related by a finite sequence of $C_{n+1}$-moves.
\end{theorem}

\begin{remark} The statement of theorem~2 remains true if we replace $C_{n+1}$-moves by \emph{special} $C_{n+1}$-moves, since these moves are equivalent (see e.g.~\cite{OTaY}).
\end{remark}
 
Hence, in our situation, the two knots $K$ and $K'$ are related by a finite sequence of special $C_{m+1}$-moves. 
By Ohyama and Tsukamoto's result~(\cite{OTs}), every single $C_{m+1}$-move changes Vassiliev invariants of order $m+1$ in a way that only depends on the cyclic permutation defined by the outer connections and a product of signs of certain crossings. By the above observation on the permutation of strands, we can imitate this pattern in our $C_{m+1}$-move from $K'$ to $\widetilde{K}'$. Repeating this procedure for every $C_{m+1}$-move 
in the sequence from $K$ to $K'$, we end up with a knot $K''$ whose Vassiliev invariants of order up to $m+1$ coincide with those of $K$. This completes the proof of theorem~1.
\end{proof}

\bigskip
\noindent
Department of Mathematics,
ETH Z\"urich, 
Switzerland

\bigskip
\noindent
\emph{sebastian.baader@math.ethz.ch}

\end{document}